\documentclass[12pt,runningheads,a4paper]{llncs}
\usepackage{amsmath,amssymb}
\usepackage{graphicx}
\usepackage{color}

\begin{document}

\title{A New Version of the M\'{e}nages Problem}
\author{Ahmad Mahmood Qureshi}
\institute{School of Mathematical Sciences,\\
GC University Lahore, Pakistan \\
sirahmad@gmail.com}
\maketitle

\newcommand{\COM}[2]{{#1\choose#2}}

\begin{abstract}
The \emph{probl\`{e}me des m\'{e}nages} (\emph{married couples
problem}) introduced by E.Lucas in 1891 is a classical problem that
asks the number of ways to arrange \emph{n} married couples around a
circular table, so that husbands and wives are in alternate places
but no couple is seated together. In this paper we present a new
version of m\'{e}nages problem that carries the constraints
consistent with several cultures.
\end{abstract}
The following problem, introduced by Lucas in 1891, is known as the
\emph{probl\`{e}me des m\'{e}nages}: in how many ways can one seat
$n$ couples at a circular table so that men and women are in
alternate places and no husband will sit on either side of his wife?
(see [1], [2]). In this paper we consider the following related problem:\\\\
Suppose there are $n_k$ \emph{k}-tuples (one husband and $k-1$
wife/wives),
$2\leq k \leq r$, $m$ single men and $w$ single women. In how many ways we can seat them around a circular table such that:\\
 1.  All members of a family sit together;\\
 2.  No man sits adjacent to a woman except his own wife/wives;\\
 3.  A husband having more than one wife will be surrounded by his wives.\\
 For ease of notation, we restrict ourselves to the case when $r=5$ (a real situation in muslim culture). In this case there
 are:\\
 $n_2$ couples (husband and wife),\\
 $n_3$ triples (husband and two wives),\\
 $n_4$ quadruples (husband and three wives),\\
 $n_5$ pentuples (husband and four wives),\\
 $m$ single men and $w$ single women.\\
\newpage
\noindent \textbf{Theorem 1.1.} \emph{The solution to this problem
for $r=5$ and $n_2$ even is}

$$(1\cdot2!)^{n_3}(2\cdot3!)^{n_4}(3\cdot4!)^{n_5} \cdot
\frac{n_2!}{\left(\frac{n_2}{2}\right)!}\cdot
\frac{\left(\frac{n_2}{2}+m-1\right)!}{\left(\frac{n_2}{2}-1\right)!}\cdot
\left(\frac{n_2}{2}+n_3+n_4+n_5+w-1\right)!$$

\noindent \textbf{Proof:}

\noindent Let us code a man by 1 and a woman by 0.
  We shall first find the number of ways to arrange the married persons around the circular table.
  Since all members of a family have to sit together, the possible coding arrangements of the
  given tuples (without distinguishing the wives) could be:\\

$\begin{tabular}{llllll}
    Couple & \indent01 & \indent10 &  &  &  \\
    Triple & \indent001 & \indent010 & \indent100 &  &  \\
    Quadruple& \indent0001 &\indent0010 &\indent 0100 &\indent 1000 &  \\
    Pentuple & \indent00001 &\indent 00010&\indent 00100 &\indent 01000 &\indent 10000 \\
  \end{tabular}$\\\\
\noindent To fulfill constraints 1 and 2, we have a unique
possibility to combine the arrangements having codes 01 and 10 from
the couples to form the block having code 0110. For the others, if
we merge consecutive zeros as one zero, then constraint 3 would only
allow blocks having codes 010 with 3, 4 or 5 members representing a single
family.\\\\
Let $i$ denote the number of 0110 blocks,
and $j$ denote the number of 010 blocks in a given arrangement. \\
Counting the number of men in all the blocks, we find
\begin{eqnarray*}
  2i + j &=& n_2 + n_3 + n_4 + n_5
\end{eqnarray*}
From the conditions of the problem it follows that the number of
solutions to this problem is different from zero if $n_2$ is even
and in this case $i=\frac{n_2}{2}$, which implies $j=n_3 + n_4 +
n_5$.\\
The number of circular arrangements of these $i+j$ blocks is
$$(i + j - 1)!\\
=\left(\frac{n_2}{2} + n_3 + n_4 + n_5 - 1\right)!$$
 Besides the circular arrangements of blocks,  there are several ways
to arrange the persons within a block.\\

\noindent \emph{Number of ways of choosing 010 blocks}:\\
For a \emph{k}-tuple ($1$ husband and $k-1$ wives ) the husband may
occupy any of the $(k-2)$ intermediate positions, while the wives
may be positioned in the block in $(k-1)!$ ways. This gives a total
of $(k-2)\cdot(k-1)!$ internal arrangements of a 010 block for a
\emph{k}-tuple.
\newpage
\noindent Accordingly, a triple may be coded as 010 in $1\cdot2!$ ways;\\
a quadruple may be coded as 010 in $2\cdot3!$ ways;\\
a pentuple may be coded as 010 in $3\cdot4!$ ways.
\\Since there are $n_3$ triples, $n_4$ quadruples and $n_5$ pentuples, the
total number of ways of choosing 010 blocks is

$$(1\cdot2!)^{n_3}(2\cdot3!)^{n_4}(3\cdot4!)^{n_5}$$

\noindent \emph{Number of ways of choosing 0110 blocks}:\\
Any 0110 block consists of two couples, so there are
$\frac{n_2}{2}=i$ such blocks.\\
The number of ways of pairing couples to generate these
$\frac{n_2}{2}$ blocks is
$$\frac{\binom{n_2}{2}\binom{n_2-2}{2}\cdots\binom{4}{2}\binom{2}{2}}{\left(\frac{n_2}{2}\right)!}
= \frac{n_2!}{\left(\frac{n_2}{2}\right)! \cdot 2^{\frac{n_2}{2}}}$$

\noindent Within each block, there are two ways of arranging the two
couples $(W_1H_1H_2W_2\,\,$ or $\,\, W_2H_2H_1W_1)$. Therefore, the
number of ways of internal arrangements in all $\frac{n_2}{2}$
blocks is $2^{\frac{n_2}{2}}$ and so the total number of ways of
choosing 0110 blocks is
$\frac{n_2!}{(\frac{n_2}{2})!}$.\\\\
Hence, the number of ways of arranging the married persons around
the circular table as per the constraints is:

$$\left(\frac{n_2}{2} + n_3+n_4 + n_5 - 1\right)! \cdot
(1\cdot2!)^{n_3}(2\cdot3!)^{n_4}(3\cdot4!)^{n_5}\cdot
\frac{n_2!}{\left(\frac{n_2}{2}\right)!}$$

\noindent Finally, the \emph{m} single men may only be placed
in-between two male persons of the \emph{i} 0110 blocks. If $r_1,
r_2,\cdots ,r_i$ men are placed in these allowed \emph{i} slots,
then
\begin{eqnarray*}
  r_1+r_2+\cdots +r_i &=& m,
\end{eqnarray*}
where $r_k \geq 0$,\,\, $1\leq k\leq i$.\\
The number of integer solutions of this equation is \\
$$\binom{i+m-1}{m}=\frac{(i+m-1)!}{m! \cdot (i-1)!}$$
\\
Any of these solutions represents a way to distribute in the $i$
slots each of the $m!$ permutations of the single men. Hence, the
number of all possible arrangements in which $m$ men, respecting the
constraints, can be seated around the table for any fixed
arrangement of married people is:
\begin{eqnarray*}
  \frac{(i+m-1)!}{m! \cdot (i-1)!}\cdot m! &=& \frac{(i+m-1)!}{(i-1)!} \\
  &=& \frac{(\frac{n_2}{2}+m-1)!}{(\frac{n_2}{2}-1)!} \\
  \end{eqnarray*}

\noindent Similarly, the \emph{w} single women may only be placed
in-between blocks. Accordingly, there are $i+j$ slots available for
them. Hence the number of arrangements in which $w$ single women,
obeying the constraints, can be seated around the table for any
fixed arrangement of married people and single men is:
\begin{eqnarray*}
  \frac{(i+j+w-1)!}{(i+j-1)!}&=& \frac{(\frac{n_2}{2}+n_3+n_4+n_5+w-1)!}{(\frac{n_2}{2}+n_3+n_4+n_5-1)!}\\
  \end{eqnarray*}

\noindent Hence the number of ways of seating all persons around the
circular table considering opposite orientations different and
respecting the given constraints with $n_2$ even is:

{\small $$(\frac{n_2}{2} + n_3 +n_4 + n_5 - 1)! \cdot
(1\cdot2!)^{n_3}(2\cdot3!)^{n_4}(3\cdot4!)^{n_5}\cdot
\frac{n_2!}{(\frac{n_2}{2})!}\cdot
\frac{(\frac{n_2}{2}+m-1)!}{(\frac{n_2}{2}-1)!}\cdot
\frac{(\frac{n_2}{2}+n_3+n_4+n_5+w-1)!}{(\frac{n_2}{2}+n_3+n_4+n_5-1)!}$$}
{\small$$=(1\cdot2!)^{n_3}(2\cdot3!)^{n_4}(3\cdot4!)^{n_5} \cdot
\frac{n_2!}{\left(\frac{n_2}{2}\right)!}\cdot
\frac{\left(\frac{n_2}{2}+m-1\right)!}{\left(\frac{n_2}{2}-1\right)!}\cdot
\left(\frac{n_2}{2}+n_3+n_4+n_5+w-1\right)! \,\,\,\,\,\,\,\,\,\,\,\, \blacksquare $$}\\

\noindent \textbf{Corollary 1.2.} \emph{The solution to this problem
for general $r$ and $n_2$ even is}\\
$$\prod_{k=3}^{r}((k-2)\cdot(k-1)!)^{n_k}\cdot\frac{n_2!}{\left(\frac{n_2}{2}\right)!}\cdot
\frac{\left(\frac{n_2}{2}+m-1\right)!}{\left(\frac{n_2}{2}-1\right)!}\cdot
\left(\sum_{k=3}^{r}{n_k}+\frac{n_2}{2}+w-1\right)!$$\\

\end{document}